\documentclass[11pt]{amsart}
   \usepackage{amsmath,amssymb}
   
   \renewcommand{\bf}{\bfseries}

   \newtheorem{theorem}{Theorem}
   \newtheorem{lemma}{Lemma}

   \renewcommand{\epsilon}{\varepsilon}

   \newcommand{\dis}{\displaystyle}
  
\begin{document}
  \title{Pointwise convergence of averages along cubes}
   \author{I. Assani}
 \begin{abstract}
   Let $(X,\mathcal{B},\mu, T)$ be a measure preserving system.
   We prove the pointwise convergence of the averages
    $$\frac{1}{N^2}\sum_{n,m= 0}^{N-1}
    f_1(T^nx)f_2(T^mx)f_3(T^{n+m}x)$$ and of similar averages
    with seven bounded functions.
  \end{abstract}
  \maketitle
  \section{Introduction}

  In \cite{Bergelson}, V. Bergelson generalized Khintchine's theorem
   \cite{Khintchine} by
  proving the $L^2$ convergence of the averages
$$\frac{1}{N^2}\sum_{n,m= 0}^{N-1}
    f_1(T^nx)f_2(T^mx)f_3(T^{n+m}x)$$ where the functions $f_i$
    are bounded measurable and $(X, \mathcal{B}, \mu, T)$ is a
    measure preserving system.
In \cite{Host-Kra}, B. Host and B. Kra extended his result by
proving the $L^2$ convergence of the following averages
  $$\frac{1}{N^3} \sum_{m,n,p = 0}^{N-1}
  f_1(T^mx)f_2(T^nx)f_3(T^{m+n}x)f_4(T^px)f_5(T^{m+p}x)f_6(T^{n+p}x)f_7(T^{m+n+p}x)$$
   They also proved that if $T$ is ergodic and all functions $f_i$ are in the
 $CL$ factor for $T$ then the averages of these
   seven functions converge a.e.. The pointwise convergence on
   such factors is a consequence of A. Leibman's result
   \cite{Leibman}

We want to show that these averages actually converge a.e. by
showing the a.e. convergence when one of the functions $f_i$
belongs to $CL^{\perp}$.

\begin{theorem}
 Let $(X,\mathcal{B},\mu, T)$ be a measure preserving system. If the functions $f_i$, $1\leq i\leq 7$, are all bounded
 then the averages
 $$ \frac{1}{N^3} \sum_{m,n,p = 0}^{N-1}
  f_1(T^mx)f_2(T^nx)f_3(T^{m+n}x)f_4(T^px)f_5(T^{m+p}x)f_6(T^{n+p}x)f_7(T^{m+n+p}x)$$
  converge a.e.
 \end{theorem}
 A corollary of our method of proof is the following result.
 \begin{theorem}
 Let $(X, \mathcal{B}, \mu, T)$ be an ergodic dynamical system.
 Then
 \begin{enumerate}
  \item Its Kronecker factor is characteristic for the pointwise
  convergence of the averages $\dis \frac{1}{N^2}\sum_{n,m= 0}^{N-1}
    f_1(T^nx)f_2(T^mx)f_3(T^{n+m}x)$
 \item Its $CL$ factor is characteristic for the pointwise
 convergence of the averages $$\frac{1}{N^3} \sum_{m,n,p = 0}^{N-1}
  f_1(T^mx)f_2(T^nx)f_3(T^{m+n}x)f_4(T^px)f_5(T^{m+p}x)f_6(T^{n+p}x)f_7(T^{m+n+p}x).$$
 \end{enumerate}
 \end{theorem}
 The notion of characteristic factor is originally due to H.
 Furstenberg. It is explicitly stated in \cite{Furstenberg-Weiss}.
 In the weakly mixing case we have the following result.
\begin{theorem}
 Let $(X, \mathcal{B}, \mu, T)$ be a weakly mixing dynamical
 system. The averages $\dis M_N(f_1, f_2,..., f_{2^k -1})$
of $2^{k}-1 $ bounded functions $f_i$ converge a.e. to $\dis
\prod_{i=1}^{2^{k}-1}\int f_i d\mu$ for all $k\geq 1.$
\end{theorem}
\section{Proofs}
 In the subsequent inequalities the constant $C$  may change
 from one line to the other. It will depend only at time on the $L^{\infty}$ norm of the functions $f_j$.

 \subsection{ Pointwise convergence for the averages of three
 functions}

\smallskip

We start by proving the pointwise convergence of the averages
$$M_N(f_1, f_2, f_3)(x)= \frac{1}{N^2}\sum_{n,m= 0}^{N-1}
    f_1(T^nx)f_2(T^mx)f_3(T^{n+m}x)$$ for $f_i$ bounded and
    measurable functions. This will help illustrate the method.
 We assume without loss of generality that $T$ is ergodic.
 We recall Bourgain's uniform Wiener Wintner ergodic result announced
 in \cite{Bourgain}.

 \begin{lemma}
  Let $(X, \mathcal{B}, \mu, T)$ be an ergodic dynamical system
  and $f$ a function in the orthocomplement of the Kronecker
  factor. Then for a.e. x we have
  $\dis \lim_N \sup_{t}|\frac{1}{N}\sum_{n=0}^{N-1} f(T^nx) e^{2\pi
  int}|= 0.$
 \end{lemma}

 Using this lemma we can prove the following

 \begin{theorem}
  Let $(X, \mathcal{B},\mu, T)$ be a measure preserving system
  and $f_i$, $1\leq i\leq 3$ three bounded functions then the averages
$$M_N(f_1, f_2, f_3)(x)= \frac{1}{N^2}\sum_{n,m= 0}^{N-1}
    f_1(T^nx)f_2(T^mx)f_3(T^{n+m}x)$$
converge a.e.
\end{theorem}

\begin{proof}

It is enough to show this convergence for ergodic measure
preserving systems (using the ergodic decomposition). We have the
following inequalities.

\[
\begin{aligned}
&|M_N(f_1, f_2, f_3)(x)|^2 \\
&\leq
\|f_1\|_{\infty}^2\bigg(\frac{1}{N}\sum_{n=0}^{N-1}\big|\frac{1}{N}\sum_{m=0}^{N-1}f_2(T^mx)f_3(T^{n+m}x)\big|^2\bigg)
\\
&\leq \|f_1\|_{\infty}^2 \frac{1}{N}\sum_{n=0}^{N-1}\bigg|\int
\big(\sum_{m=0}^{N-1}f_2(T^mx)e^{-2\pi
imt}\big)\big(\frac{1}{N}\sum_{m'=0}^{2(N-1)}f_3(T^{m'}x)e^{2\pi
im't}\big).
e^{2\pi int} dt \bigg|^2 \\
&\leq \|f_1\|_{\infty}^2\frac{1}{N}\int \bigg|
\sum_{m=0}^{N-1}f_2(T^mx)e^{-2\pi
imt}\bigg|^2\bigg|\frac{1}{N}\sum_{m'=0}^{2(N-1)}f_3(T^{m'}x)e^{2\pi
im't}\bigg|^2dt \\
& \leq
\frac{C}{N}\sup_t\bigg|\frac{1}{N}\sum_{m'=0}^{N-1}f_3(T^{m'}x)e^{2\pi
im't}\bigg|^2\int \big|\sum_{m=0}^{N-1}f_2(T^{m}x)e^{-2\pi
imt}|^2dt \\
& \leq
C\sup_t\bigg|\frac{1}{N}\sum_{m'=0}^{N-1}f_3(T^{m'}x)e^{2\pi
im't}\bigg|^2\frac{1}{N} N\|f_2\|_{\infty}^2
\end{aligned}
\]

With the help of lemma 1 we can conclude that for $f_3$ in the
orthocomplement of the Kronecker factor the averages
$M_N(f_1,f_2,f_3)$ converge a.e. to zero.

If $f_3$ is one of the eigenfunctions for $T$ with eigenvalue
$e^{2\pi i\theta}$ then
$$M_N(f_1,f_2, f_3) = f_3\big(\frac{1}{N}\sum_{n=0}^{N-1}f_1(T^nx)e^{2\pi
in\theta}\big)\big(\frac{1}{N}\sum_{m=0}^{N-1}f_2(T^mx) e^{2\pi
im\theta}\big).$$ The convergence in this case follows from
Birkhoff 's theorem applied to the product of T and the rotation
$\theta$.
 The convergence for a general function $f_3$ in the Kronecker
 factor follows now by linearity and approximation.

\end{proof}
\noindent{\bf Remarks 1}

\vskip1ex
\begin{itemize}
\item
  The proof of theorem 4 shows that
  if $f_1$ and $f_2$ are bounded functions and $P_{\mathcal{K}}$
  denotes the projection onto the Kronecker factor of $T$ then
  \begin{equation}
  \limsup_N\frac{1}{N}\sum_{n=0}^{N-1}\big|\frac{1}{N}\sum_{m=0}^{N-1}f_1(T^mx)f_2(T^{m+n}x)\big|^2
  =
  \limsup_N\frac{1}{N}\sum_{n=0}^{N-1}\big|\frac{1}{N}\sum_{m=0}^{N-1}P_{\mathcal{K}}(f_1)(T^mx)P_{\mathcal{K}}(f_2)(T^{m+n}x)\big|^2
  \end{equation}
\item
  The proof of theorem 4 actually shows that
  \begin{equation}
  \bigg(\frac{1}{N}\sum_{n=0}^{N-1}\big|\frac{1}{N}\sum_{m=0}^{N-1}f_2(T^mx)f_3(T^{n+m}x)\big|^2\bigg)\leq
  C\sup_t\bigg|\frac{1}{N}\sum_{m'=0}^{N-1}f_3(T^{m'}x)e^{2\pi
im't}\bigg|^2\|f_2\|_{\infty}^2.
  \end{equation}
A similar estimate can be obtained with $\dis
\sup_t\bigg|\frac{1}{N}\sum_{m'=0}^{N-1}f_2(T^{m'}x)e^{2\pi
im't}\bigg|^2$ if we focus instead on the function $f_2$.
\end{itemize}
\subsection{Pointwise convergence for the averages of seven
functions}

\smallskip

  As $T$ is ergodic there exists in $\mathcal{K}$ an orthonormal
 basis of eigenfunctions $g_j$ with modulus 1 corresponding to the eigenvalue $e^{2\pi i\theta_j}$
 so that any function $G\in
 \mathcal{K}$ can be written as
 \begin{equation}
 G= \sum_{j=1}^{\infty} \big(\int G.\overline{g_j} d\mu\big) g_j.
 \end{equation}
  In \cite{Host-Kra1} it is shown that the CL factor is characteristic for the convergence in $L^2$ norm of the averages of seven functions. Functions in this factor
  are characterized by the seminorm $|\|.|\|_3$ such that
  \begin{equation}
  \||f|\|_3^8 = \lim_{H}\frac{1}{H} \sum_{h=0}^{H-1}\||f.f\circ
  T^h|\|_2^4
  \end{equation}
   where
  \begin{equation}
  \||f|\|_2^4 = \lim_{H}\frac{1}{H}\sum_{h=0}^{H-1} \big|\int f.f(T^h)d\mu\big|^2.
  \end{equation}
     A function $f\in CL^{\perp}$ if and only
  $\||f|\|_3 =0$.

  \begin{lemma}
  Let $(X, \mathcal{B}, \mu, T)$ be an ergodic dynamical system
  and $f\in L^{\infty}(\mu)$ then for all $H$ positive integer
  we have
  $$\limsup_N \sup_t \big|\frac{1}{N}\sum_{n=0}^{N-1}f(T^nx)e^{2\pi
  int}\big|^2 \leq C\bigg( \frac{1}{H} + \frac{1}{H}\sum_{h=1}^H
  \big|\int f.\overline{f\circ T^h}d\mu\big|\bigg)
  $$
  In particular we have
  \begin{equation}
  \limsup_N \sup_t \big|\frac{1}{N}\sum_{n=0}^{N-1}f(T^nx)e^{2\pi
  int}\big|^2 \leq C\||f|\|_2^2.
  \end{equation}
  \end{lemma}
 \begin{proof}
 Without loss of generality we can assume that the function $f$
 takes only real values.
 We apply van der Corput's inequality
 (\cite{Kuipers-Niederreiter}) . For $H<N$ we get
 $$\sup_t\big|\frac{1}{N}\sum_{n=0}^{N-1} f(T^nx)e^{2\pi
 int}\big|^2 \leq C\bigg(\frac{1}{H} + \frac{1}{H}\sum_{h=1}^H
 \bigg|\frac{1}{N}\sum_{n=0}^{N-h} f(T^nx)f(T^{n+h}x)\bigg|\bigg)
 $$
  Birkhoff's pointwise ergodic theorem allows us to obtain the
 first part of the lemma.
 \smallskip
 For the second part we can use Cauchy Schwartz inequality to
 write that
$$\frac{1}{H}\sum_{h=1}^H
 \big|\int f.f\circ T^h d\mu \big|
\leq \bigg(\frac{1}{H}\sum_{h=1}^H\big|\int f.f\circ T^h d\mu
\big|^2\bigg)^{1/2}.
 $$
Now using the definition of $\||f|\|_2$, (see (5)), we can end the
proof of this lemma.
 \end{proof}

 The lemma that replaces the uniform Wiener Wintner
ergodic theorem in the case of the averages of seven functions is
the following.

\begin{lemma}
If $f_1$ or $f_2$ is in $CL^{\perp}$ then for a.e. $x$
\begin{equation}
\lim_N\frac{1}{N}\sum_{n=0}^{N-1}\sup_t\bigg|\frac{1}{N}\sum_{m=0}^{N-1}f_1(T^mx)f_2(T^{n+m}x)e^{2\pi
imt}\bigg|^2 = 0
\end{equation}
\end{lemma}
\begin{proof}

We can assume without loss of generalities that the functions are
uniformly bounded by one. We use again van der Corput's
inequality, \cite{Kuipers-Niederreiter}. For $(H +1)^2< N$ we get
\[
\begin{aligned}
&\sup_t\bigg|\frac{1}{N}\sum_{m=0}^{N-1}f_1(T^mx)f_2(T^{n+m}x)e^{2\pi
imt}\bigg|^2  \\
&\leq \frac{C}{H} +
\frac{C}{H}\sum_{h=1}^H\bigg|\frac{1}{N}\sum_{m=0}^{N-h-1}f_1(T^mx)f_2(T^{m+n}x)\overline{f_1(T^{m+h}x)}
\overline{f_2(T^{m+n+h}x)}\bigg|
\end{aligned}
\]
So recalling that the constant C may change from one line to
another but remains an absolute constant we have,
\[
\begin{aligned}
&\frac{1}{N}\sum_{n=0}^{N-1}\sup_t\bigg|\frac{1}{N}\sum_{m=0}^{N-1}f_1(T^mx)f_2(T^{n+m}x)e^{2\pi
imt}\bigg|^2 \\
&\leq \frac{C}{H} +
\frac{C}{H}\sum_{h=1}^H\frac{1}{N}\sum_{n=0}^{N-1}\bigg|\frac{1}{N}\sum_{m=0}^{N-h-1}f_1(T^mx)f_2(T^{m+n}x)\overline{f_1(T^{m+h}x)}
\overline{f_2(T^{m+n+h}x)}\bigg| \\
&\leq \frac{C}{H} +
\frac{C}{H}\sum_{h=1}^H\frac{1}{N}\sum_{n=0}^{N-1}\bigg|\frac{1}{N}\sum_{m=0}^{N-1}f_1(T^mx)f_2(T^{m+n}x)\overline{f_1(T^{m+h}x)}
\overline{f_2(T^{m+n+h}x)}\\
&-\sum_{m=N-h}^{N-1}f_1(T^mx)f_2(T^{m+n}x)\overline{f_1(T^{m+h}x)}
\overline{f_2(T^{m+n+h}x)}\bigg| \\
&\leq \frac{C}{H} +
\frac{C}{H}\sum_{h=1}^H\frac{1}{N}\sum_{n=0}^{N-1}\bigg|\frac{1}{N}\sum_{m=0}^{N-1}f_1(T^mx)f_2(T^{m+n}x)\overline{f_1(T^{m+h}x)}
\overline{f_2(T^{m+n+h}x)}\bigg| +
\frac{C}{H}\sum_{h=1}^H\frac{1}{N}\sum_{n=0}^{N-1}\frac{h}{N}\\
&\leq \frac{C}{H} +
\frac{C}{H}\sum_{h=1}^H\frac{1}{N}\sum_{n=0}^{N-1}\bigg|\frac{1}{N}\sum_{m=0}^{N-1}f_1(T^mx)f_2(T^{m+n}x)\overline{f_1(T^{m+h}x)}
\overline{f_2(T^{m+n+h}x)}\bigg|.
\end{aligned}
\]
 Thus using the inequality ( or Cauchy Schwartz's inequality)
 \begin{equation}
|\frac{1}{P}\sum_{p=1}^P u_p|\leq (\frac{1}{P}\sum_{p=1}^P
|u_p|^2)^{1/2}
\end{equation}
 we obtain
 \[
 \begin{aligned}
&\frac{1}{N}\sum_{n=0}^{N-1}\sup_t\bigg|\frac{1}{N}\sum_{m=0}^{N-1}f_1(T^mx)f_2(T^{n+m}x)e^{2\pi
imt}\bigg|^2 \\
&\leq \frac{C}{H} + \bigg(\frac{C}{H}\sum_{h=1}^H
\bigg(\frac{1}{N}\sum_{n=0}^{N-1}\bigg|\frac{1}{N}\sum_{m=0}^{N-1}f_1(T^mx)f_2(T^{m+n}x)\overline{f_1(T^{m+h}x)}
\overline{f_2(T^{m+n+h}x)}\bigg|^2\bigg)\bigg)^{1/2}
\end{aligned}
\]
Finally by applying the inequality (2) made after the Remarks 1 to
the function $\dis f_1.\overline{f_1\circ T^h}$ we get
\[
\begin{aligned}
&\frac{1}{N}\sum_{n=0}^{N-1}\sup_t\bigg|\frac{1}{N}\sum_{m=0}^{N-1}f_1(T^mx)f_2(T^{n+m}x)e^{2\pi
imt}\bigg|^2 \\
&\leq \frac{C}{H} + \bigg(\frac{C}{H}\sum_{h=1}^H
\bigg(\sup_t\bigg|\frac{1}{N}\sum_{m'=0}^{N-1}(f_1.\overline{f_1\circ
T^h})(T^{m'}x)e^{2\pi im't}\bigg|^2\bigg)^{1/2}
\end{aligned}
\]

Now by using Lemma 2 and the inequality $\dis
\frac{1}{H}\sum_{h=1}^H |u_h|^2\leq \big(\frac{1}{H}\sum_{h=1}^H
|u_h|^4\big)^{1/2}$ we obtain
\[
\begin{aligned}
&\limsup_N\frac{1}{N}\sum_{n=0}^{N-1}\sup_t\bigg|\frac{1}{N}\sum_{m=0}^{N-1}f_1(T^mx)f_2(T^{n+m}x)e^{2\pi
imt}\bigg|^2 \\
&\leq \frac{C}{H} + \bigg(\frac{C}{H}\sum_{h=1}^H
\limsup_N\sup_t\bigg|\frac{1}{N}\sum_{m'=0}^{N-1}(f_1.\overline{f_1\circ
T^h})(T^{m'}x)e^{2\pi im't}\bigg|^2\bigg)^{1/2}\\
&\leq \frac{C}{H} + \bigg(\frac{C}{H}\sum_{h=1}^H
\||f_1.\overline{f_1\circ T^h}|\|_2^2\bigg)^{1/2}\\
&\leq \frac{C}{H} + \bigg(\frac{C}{H}\sum_{h=1}^H
\||f_1.\overline{f_1\circ T^h}|\|_2^2\bigg)^{1/2}\\
&\leq \frac{C}{H} + \bigg(\frac{C}{H}\sum_{h=1}^H
\||f_1.\overline{f_1\circ T^h}|\|_2^4\bigg)^{1/4}
\end{aligned}
\]
Taking now the limit when $H$ tends to $\infty$ we get the
following estimate
\begin{equation}
\limsup_N\frac{1}{N}\sum_{n=0}^{N-1}\sup_t\bigg|\frac{1}{N}\sum_{m=0}^{N-1}f_1(T^mx)f_2(T^{n+m}x)e^{2\pi
imt}\bigg|^2\leq C\||f_1|\|_3^2
\end{equation}

 Thus if we assume that $f_1\in CL^{\perp}$ then $\||f_1|\|_3 =0$ and we
 obtain the equation (7). We have the same conclusion if one assumes
 that $f_2\in CL^{\perp}$.
\end{proof}
 Using Lemma 3 we can now give a proof of theorem 1.

 \begin{proof}{\bf Theorem 1}
\[
\begin{aligned}
&|M_N(f_1,f_2,...,f_7)|^2 \\
&=\bigg|
\frac{1}{N^3}\sum_{p=0}^{N-1}f_1(T^px)\sum_{n=0}^{N-1}f_2(T^nx)f_3(T^{p+n}x)\big(\sum_{m=0}^{N-1}f_4(T^mx)f_5(T^{n+m}x)f_6(T^{p+m}x)
f_7(T^{n+m+p}x)\big)\bigg|^2 \\
&\leq
\frac{1}{N^2}\sum_{p=0}^{N-1}\sum_{n=0}^{N-1}\|f_1\|_{\infty}^2\|f_2\|_{\infty}^2\|f_3\|_{\infty}^2\bigg|\frac{1}{N}
\sum_{m=0}^{N-1}f_4(T^mx)f_5(T^{n+m}x)f_6(T^{p+m}x)f_7(T^{p+n+m}x)\bigg|^2
\\
&= \frac{1}{N^2}\prod_{i=1}^{3}\|f_i\|_{\infty}^2. \\
&\sum_{n=0}^{N-1}\sum_{p=0}^{N-1}\bigg|\int \big(
\sum_{m=0}^{(N-1)} f_4(T^mx)f_5(T^{n+m}x)e^{-2\pi
imt}\big)\big(\frac{1}{N}\sum_{m'=0}^{2(N-1)}f_6(T^{m'}x)f_7(T^{n+m'}x)e^{2\pi
im't}\big). e^{2\pi ipt}dt\bigg|^2 \\
&\leq
\frac{1}{N^2}\prod_{i=1}^3\|f_i\|_{\infty}^2\sum_{n=0}^{N-1}\int
\bigg| \sum_{m=0}^{N-1} f_4(T^mx)f_5(T^{n+m}x)e^{-2\pi
imt}\big)\big(\frac{1}{N}\sum_{m'=0}^{2(N-1)}f_6(T^{m'}x)f_7(T^{n+m'}x)e^{2\pi
im't}\big)\bigg|^2 dt \\
&\leq
\frac{C}{N^2}\prod_{i=1}^3\|f_i\|_{\infty}^2\sum_{n=0}^{N-1}\sup_t
\bigg|\frac{1}{N}\sum_{m'=0}^{N-1}f_6(T^{m'}x)f_7(T^{n+m'}x)e^{2\pi
im't}\bigg|^2 N\prod_{j=4}^{5}\|f_j\|_{\infty}^2 \\
&= C\prod_{i=1}^{5}\|f_i\|_{\infty}^2
\frac{1}{N}\sum_{n=0}^{N-1}\sup_t
\bigg|\frac{1}{N}\sum_{m'=0}^{N-1}f_6(T^{m'}x)f_7(T^{n+m'}x)e^{2\pi
im't}\bigg|^2
\end{aligned}
\]
With the help of lemma 3 one can conclude that if $f_6$ or $f_7$
belong to $CL^{\perp}$ then the averages of these seven functions
converge to zero. By using the symmetry of the sum of the averages
with respect to $n$, $m$ and $p$ one can see that the averages
will converge to zero if one of the functions $f_i\in CL^{\perp},$
$1\leq i\leq 7$.

 \end{proof}
\noindent{\bf Remarks 2}
\begin{itemize}
 \item
 The last steps of the proof of theorem 1 show that for bounded functions $f_i$,
 $4\leq i \leq 7$ if we denote by $P_{CL}(f_i)$ their projection onto the $CL$ factor then we have
 \begin{equation}
 \begin{aligned}
 &\limsup_N \frac{1}{N^2}\sum_{n,p =0}^{N-1}\bigg|\frac{1}{N}
\sum_{m=0}^{N-1}f_4(T^mx)f_5(T^{n+m}x)f_6(T^{p+m}x)f_7(T^{p+n+m}x)\bigg|^2\\
&= \limsup_N\frac{1}{N^2}\sum_{n,p =0}^{N-1}\bigg|\frac{1}{N}
\sum_{m=0}^{N-1}P_{CL}(f_4)(T^mx)P_{CL}(f_5)(T^{n+m}x)P_{CL}(f_6)(T^{p+m}x)P_{CL}(f_7)(T^{p+n+m}x)\bigg|^2.
\end{aligned}
\end{equation}
\item
The proof of lemma 3 gives the following estimate
\begin{equation}
\limsup_N\frac{1}{N}\sum_{n=0}^{N-1}\sup_t\bigg|\frac{1}{N}\sum_{m=0}^{N-1}f_1(T^mx)f_2(T^{n+m}x)e^{2\pi
imt}\bigg|^2\leq C Min[\||f_1|\|_3^2, \||f_2|\|_3^2].
\end{equation}
\end{itemize}

 \subsection{Proof of Theorem 2.}
 \vskip1ex

  The proof is a consequence of the path used in establishing theorem
  1. We have shown that if one of the functions $f_i\in
  CL^{\perp}$, $1\leq i\leq 7$,  then the averages converge
  pointwise to zero. This shows that the $CL$ factor is
  characteristic for the pointwise convergence. For the averages
  of three functions the Kronecker factor is characterisitc for
  the pointwise convergence for the same reason.

\subsection{Proof of Theorem 3.}
\vskip1ex

  We list some properties and some notations. They may seem a bit
  complicated at first reading. So the reader may wish to first translate all these properties
  to the case of 15 functions.
\begin{enumerate}
 \item For each $k\geq 4$ we denote by
 $$M_N(f_1, f_2,..., f_{2^k -1})(x)$$ the averages of
 $2^k-1$ bounded functions . We number the functions $f_j$ so that
 those with $2^{k-1}\leq j\leq 2^k -1$  are depending of
the index $i_k$.  For instance in the sum of 7 functions, the
 functions are $f_j$, $4\leq j\leq 7$ and they appear in the sum
 $\dis \sum_{m=0}^{N-1}f_4(T^mx)f_5(T^{n+m}x)f_6(T^{p+m}x)f_7(T^{p+n+m}x).$
 In the case of 15 functions if we denote by $p,n,k,m $ the indices $i_1,i_2, i_3, i_4$ then
 they appear in the sum
 $$ \sum_{m=0}^{N-1}f_8(T^mx)f_9(T^{n+m}x)...f_{15}(T^{p+n+k+m}x)$$
 We denote by  $S_{N,(i_1,i_2,...,i_k)}(f_{2^{k-1}},..., f_{2^k -1})(x)$ these
 terms depending on $i_k$. We can observe that each term
   $S_{N, (i_1,i_2,...,i_k)}(f_{2^{k-1}},..., f_{2^k -1})(x)$ is
 the product of two groups of $2^{k-2}$ functions,$$
 A_{N, (i_2,...,i_{k-1},{i_k})}(f_{2^{k-1}}, f_{2^{k-1}+1},...,f_{3.2^{k-2}})(x)$$ and
 $$ B_{N, (i_1,i_2,...,i_k)}(f_{3.2^{k-2}+1},...,f_{2^k -1})(x)$$ such that the
 powers of T associated with each function in the second group are
 exactly those associated with the functions in the first group
 shifted by the index $i_{1}$.  Similar decompositions can be
 obtained if one focus on shifted blocks by another index. One can
 observe that we could write
 \begin{equation}
 B_{N, (i_1,i_2,...,i_k)}(f_{3.2^{k-2}+1},...,f_{2^k -1})(x)
 =A_{N, (i_2,...,i_{k-1},{i_k})}(f_{3.2^{k-2}+1},...,f_{2^k -1})(T^{i_1}x)
 \end{equation}
 \vskip1ex

 The interest in those terms in the numerator of $\dis M_N(f_1,
f_2,..., f_{2^k -1})(x)$ rests also in the following
\begin{equation}
\begin{aligned}
&|M_N(f_1, f_2,..., f_{2^k -1})(x)|^2 \\
&\leq \prod_{j=1}^{2^{k-1}-1}\|f_j\|_{\infty}^2
\frac{1}{N^{k-1}}\sum_{i_1,...,
 i_{k-1}=0}^{N-1}\big|\frac{1}{N}\sum_{i_k=0}^{N-1}S_{N,(i_1,i_2,...,i_k)}(f_{2^{k-1}},..., f_{2^k
-1})(x)\big|^2.
\end{aligned}
\end{equation}

 \item
 When $T$ is weakly mixing the Kronecker and $CL$ factors are
 trivial. Thus we have $P_{K}f_i= P_{CL}(f_i) = \int f_i d\mu$.
 \vskip1ex

 We want to prove theorem 3 by induction on $k$. We formulate our induction
 assumption.
\end{enumerate}

\noindent{\bf Induction Assumption}

\vskip1ex

 We assume that the following properties hold for all bounded functions
$f_j$, $1\leq j \leq k-1$.
\begin{enumerate}
\item
\[
\begin{aligned}
&\limsup_N\frac{1}{N^{k-2}}\sum_{i_1,...,i_{k-2}=0}^{N-1}\big|\frac{1}{N}\sum_{i_{k-1}=0}^{N-1}S_{N,(i_1,i_2,...,i_{k-1})}(f_{2^{k-2}},...,
f_{2^{k-1}-1})(x)\big|^2\\
&=
\limsup_N\frac{1}{N^{k-2}}\sum_{i_1,...,i_{k-2}=0}^{N-1}\prod_{j=2^{k-2}}^{2^{k-1}-1}\big|\int
f_jd\mu\big|^2 \\
&= \prod_{j=2^{k-2}}^{2^{k-1}-1}\big|\int f_jd\mu\big|^2
\end{aligned}
\]
 (Compare these equalities to the equations (1) and (10 ) in the remarks after the proofs
 for three terms and seven terms).
 \item The averages of $2^{k-1}-1$ bounded functions converge a.e.
 to the product of the integrals of these functions.
\end{enumerate}
We want to show that the same assumptions hold then for $k$.
 We can assume that all functions are real valued.
First we want to establish the following lemma
\begin{lemma}
If one of the  $2^{k-2}$ functions $f_j$, $3.2^{k-2}+1\leq j\leq
2^{k}-1$ has zero integral then
\begin{equation}
\lim_N\frac{1}{N^{k-2}}\sum_{i_1,..., i_{k-2}=
0}^{N-1}\sup_{t}\bigg|\frac{1}{N}\sum_{i_k=0}^{N-1}A_{N,
(i_1,i_2,...,i_{k-2},i_k)}(f_{3.2^{k-2}+1},...,f_{2^k
-1})(x)e^{2\pi ii_kt}\bigg|^2 =0
\end{equation}
\end{lemma}
\begin{proof}
 As previously we apply Van der Corput lemma to each term
 \[
 \sup_{t}\bigg|\frac{1}{N}\sum_{i_k=0}^{N-1}A_{N, (i_1,i_2,...,i_{k-2},i_k)}(f_{3.2^{k-2}+1},...,f_{2^k
-1})(x)e^{2\pi ii_kt}\bigg|^2
\]
We have then for each $H<N$
\[
\begin{aligned}
&\frac{1}{N^{k-2}}\sum_{i_1,..., i_{k-2}=
0}^{N-1}\sup_{t}\bigg|\frac{1}{N}\sum_{i_k=0}^{N-1}A_{N,(i_1,i_2,...,i_{k-2},i_k)}(f_{3.2^{k-2}+1},...,f_{2^k
-1})(x)e^{2\pi ii_kt}\bigg|^2 \\
& \leq \frac{1}{N^{k-2}}\sum_{i_1,..., i_{k-2}= 0}^{N-1}
C.\bigg(\frac{1}{H} + \frac{1}{H}\sum_{h=1}^H \\
&\bigg|\frac{1}{N}\sum_{i_k=1}^{N-h-1}A_{N,(i_1,i_2,...,i_{k-2},i_k)}(f_{3.2^{k-2}+1}.f_{3.2^{k-2}+1}\circ
T^h,...,f_{2^k -1}.f_{2^k -1}\circ T^h)(x)\bigg|\bigg)\\
&\leq C.\bigg(\frac{1}{H} + \frac{1}{H}\sum_{h=1}^H
\frac{1}{N^{k-2}} \\
&\sum_{i_1,..., i_{k-2}=
0}^{N-1}\bigg|\frac{1}{N}\sum_{i_k=1}^{N-h-1}A_{N,(i_1,i_2,...,i_{k-2},i_k)}(f_{3.2^{k-2}+1}.f_{3.2^{k-2}+1}\circ
T^h,...,f_{2^k -1}.f_{2^k -1}\circ T^h)(x)\bigg|\bigg)\\
&\leq C.\bigg(\frac{1}{H} + \frac{1}{H}\sum_{h=1}^H
\frac{1}{N^{k-2}}\sum_{i_1,..., i_{k-2}=
0}^{N-1}\\
&\bigg|\frac{1}{N}\sum_{i_k=1}^{N-h-1}A_{N,(i_1,i_2,...,i_{k-2},i_k)}(f_{3.2^{k-2}+1}.f_{3.2^{k-2}+1}\circ
T^h,...,f_{2^k -1}.f_{2^k -1}\circ T^h)(x)\bigg|\bigg)
\end{aligned}
\]

Then we estimate
$$\frac{1}{H}\sum_{h=1}^H
\limsup_N \frac{1}{N^{k-2}}\sum_{i_1,..., i_{k-2}=
0}^{N-1}\bigg|\frac{1}{N}\sum_{i_k=1}^{N-h-1}A_{N,(i_1,i_2,...,i_{k-2},i_k)}(f_{3.2^{k-2}+1}.f_{3.2^{k-2}+1}\circ
T^h,...,f_{2^k -1}.f_{2^k -1}\circ T^h)(x)\bigg|
$$
which by the equation (8) (in the proof of lemma 3) is less than
\[\begin{aligned}
&\frac{1}{H}\sum_{h=1}^H\limsup_N \\
&\bigg(\frac{1}{N^{k-2}}\sum_{i_1,..., i_{k-2}=
0}^{N-1}\bigg|\frac{1}{N}\sum_{i_k=1}^{N-h-1}A_{N,(i_1,i_2,...,i_{k-2},i_k)}(f_{3.2^{k-2}+1}.f_{3.2^{k-2}+1}\circ
T^h,...,f_{2^k -1}.f_{2^k -1}\circ T^h)(x)\bigg|^2\bigg)^{1/2}
\end{aligned}
\]
Now using the first induction assumption we conclude that
\[
\begin{aligned}
& \limsup_N \bigg(\frac{1}{N^{k-2}}\sum_{i_1,..., i_{k-2}=
0}^{N-1}\\
&\bigg|\frac{1}{N}\sum_{i_k=1}^{N-h-1}A_{N,(i_1,i_2,...,i_{k-2},i_k)}(f_{3.2^{k-2}+1}.f_{3.2^{k-2}+1}\circ
T^h,...,f_{2^k -1}.f_{2^k -1}\circ T^h)(x)\bigg|^2\bigg)^{1/2}\\
&= \big(\prod_{j= 2^{k-2}}^{2^{k-1}-1} \big|\int f_j.f_j\circ
T^h|^2\big)^{1/2}
\end{aligned}
\]
 As one of the functions $f_j$ let us say $g= f_{j_0}$ has integral zero and $T$ is weakly
 mixing then the spectral measure $\sigma_g$ is continuous . Thus
 we have
  $$\lim_H \frac{1}{H}\sum_{h=1}^H \big|\int g.g\circ T^h
  d\mu\big|^2 =0$$
  As the functions are bounded
  \[
 \begin{aligned}
&\frac{1}{H}\sum_{h=1}^H\limsup_N \\
&\bigg(\frac{1}{N^{k-2}}\sum_{i_1,..., i_{k-2}=
0}^{N-1}\bigg|\frac{1}{N}\sum_{i_k=1}^{N-h-1}A_{N,(i_1,i_2,...,i_{k-2},i_k)}(f_{3.2^{k-2}+1}.f_{3.2^{k-2}+1}\circ
T^h,...,f_{2^k -1}.f_{2^k -1}\circ T^h)(x)\bigg|^2\bigg)^{1/2}\\
&\leq C. \frac{1}{H}\sum_{h=1}^h \big|\int g.g\circ T^h d\mu\big|
\end{aligned}
\]
 Taking now the limit with $H$ we obtain a proof of the lemma.
\end{proof}
\noindent{\bf Remark 3} In the case of the averages of 15
functions the equation (14) in lemma 4 is
$$\lim_N\frac{1}{N^2}\sum_{p=0}^{N-1}\sum_{n=0}^{N-1}\sup_t\bigg|\frac{1}{N}
\sum_{m=0}^{N-1}f_4(T^mx)f_5(T^{n+m}x)f_6(T^{p+m}x)f_7(T^{p+n+m}x)e^{2\pi
imt}\bigg|^2=0$$

\noindent{\bf End of the proof of theorem 3}

\vskip1ex
 We just need to prove the induction at step $l=k$.
 We consider then the averages of $2^{k}-1$ functions $f_j$ and we
 use the previous observations to write

 \[
 \begin{aligned}
&|M_N(f_1, f_2,..., f_{2^k -1})(x)|^2\\
& \leq \prod_{j=1}^{2^{k-1}-1}\|f_j\|_{\infty}^2
\frac{1}{N^{k-1}}\sum_{i_1,...,
 i_{k-1}=0}^{N-1}\big|\frac{1}{N}\sum_{i_k=0}^{N-1}S_{N,(i_1,i_2,...,i_k)}(f_{2^{k-1}},..., f_{2^k
-1})(x)\big|^2
\end{aligned}
\]
 Using the equation (14) we can write
 \[
 \begin{aligned}
&\big|\frac{1}{N}\sum_{i_k=0}^{N-1}S_{N,(i_1,i_2,...,i_k)}(f_{2^{k-1}},...,
f_{2^k -1})(x)\big|^2 \\
& =\bigg|\int \big(\frac{1}{N}\sum_{i_k=0}^{N-1}A_{N,
(i_2,...,i_{k-1},{i_k})}(f_{2^{k-1}},
f_{2^{k-1}+1},...,f_{3.2^{k-2}})(x)e^{-2\pi
i_kt}\big) \\
&\big(\sum_{i_k^{'}=0}^{N-1}A_{N, (i_2,...,i_{k-1},
i_k^{'})}(f_{3.2^{k-2}+1},...,f_{2^k -1})(x)e^{2\pi
i_k^{'}t}\big).e^{2\pi ii_1t}dt\bigg|^2
\end{aligned}
\]
Hence we have
\[
 \begin{aligned}
&|M_N(f_1, f_2,..., f_{2^k -1})(x)|^2 \\
& \leq \prod_{j=1}^{2^{k-1}-1}\|f_j\|_{\infty}^2
\frac{1}{N^{k-1}}\sum_{i_1,...,
 i_{k-1}=0}^{N-1}\big|\frac{1}{N}\sum_{i_k=0}^{N-1}S_{N,(i_1,i_2,...,i_k)}(f_{2^{k-1}},..., f_{2^k
-1})(x)\big|^2 \\
&\leq \prod_{j=1}^{2^{k-1}-1}\|f_j\|_{\infty}^2
\frac{1}{N^{k-1}}\sum_{i_2,...,
 i_{k-1}=0}^{N-1}\int \bigg|\big(\sum_{i_k=0}^{N-1}A_{N,
(i_2,...,i_{k-1},{i_k})}(f_{2^{k-1}},
f_{2^{k-1}+1},...,f_{3.2^{k-2}})(x)e^{-2\pi i_kt}\big) \\
&\big(\frac{1}{N}\sum_{i_k^{'}=0}^{N-1}A_{N, (i_2,...,i_{k-1},
i_k^{'})}(f_{3.2^{k-2}+1},...,f_{2^k -1})(x)e^{2\pi
i_k^{'}t}\big)\bigg|^2 dt\\
&\leq C\frac{1}{N^{k-2}}\sum_{i_2,...,
 i_{k-1}=0}^{N-1}\sup_t\big|\frac{1}{N}\sum_{i_k^{'}=0}^{N-1}A_{N, (i_2,...,i_{k-1},
i_k^{'})}(f_{3.2^{k-2}+1},...,f_{2^k -1})(x)e^{2\pi
i_k^{'}t}\big|^2
 \end{aligned}
 \]
By using  Lemma 4 one can conclude that the averages $\dis
M_N(f_1, f_2,..., f_{2^k -1})(x)$ converge a.e to zero when one of
the functions $f_j$ has a zero integral. (using the symmetry on
the indices). From this one derives that the averages of $2^k -1$
bounded functions converge to the product of the integral of the
functions. This is part (2) of the induction assumption at level
$k$. To end the proof of the theorem we just need to observe that
the proof given for $l=k$ proves also the first assumption for
$k$.

 \noindent{\bf Remark 4}

\vskip1ex
 If one considers instead the averages
$$\frac{1}{(N-M)^2}\sum_{n,m =M}^{N}f_1(T^nx)
f_2(T^mx)f_3(T^{n+m}x)$$ where $(N-M)$ tends to $\infty$ then we
do not have a.e. convergence in general while as shown in
\cite{Bergelson} and \cite{Host-Kra} we do have convergence in
$L^2$ norm. For instance it is shown in \cite{schwartz} that for
$\beta\geq 3$ the averages
$$\frac{1}{N^{\beta-1}}\sum_{n=N^{\beta}}^{(N+1)^{\beta}}f(T^nx)$$
do not converge a.e. even if f is the characteristic function of a
set of positive measure. So in this case the Kronecker factor is
characteristic for the $L^2$ norm but not for the pointwise
convergence.

\end{document}